\documentclass[11pt, oneside, 
a4paper]{article}
\usepackage[ansinew]{inputenc}
\usepackage[spanish]{babel}
\usepackage{latexsym}
\usepackage{amsbsy}
\usepackage{amscd}
\usepackage{amsfonts}
\usepackage{amsmath}
\usepackage{amsopn}
\usepackage{amssymb}
\usepackage{amstext}
\usepackage{amsxtra}
\usepackage{url}

\usepackage{srcltx}
\usepackage{graphicx}
\newcommand{\CC}{\mathcal{C}}
\newcommand{\DD}{\mathcal{D}}
\newcommand{\EE}{\mathcal{E}}
\newcommand{\FF}{\mathcal{F}}

\newcommand{\LL}{\mathcal{L}}

\newcommand{\PP}{\mathcal{P}}

\newcommand{\TT}{\mathcal{T}}

\begin{document}
\date{} 
\title{C\'{o}mo obtener curvas con formas predeterminadas a partir de circunferencias}
\author{M.J. de la
Puente\\
Dpto. de Algebra UCM\\
\texttt{mpuente@mat.ucm.es}\\
}

\maketitle

\begin{abstract}
We produce several algebraic curves, some well--known, some new, out of circles, by means of two classical (mutually reciprocal) algebraic methods: blow--down and blow--up.
\end{abstract}

\section{Introducci\'{o}n}
Una mirada atenta a nuestro alrededor nos revela   que vivimos en  un mundo  poblado de curvas.
En los fen\'{o}menos naturales aparecen curvas de distinta \'{\i}ndole: se forman \emph{circunferencias} conc\'{e}ntricas  al arrojar una piedra a una masa de agua en calma, las \'{o}rbitas planetarias son  \emph{elipses} y
las caracolas son \emph{espirales}. El arco iris, de delicados colores,  es otro ejemplo fascinante de arco de curva. Vemos  curvas en la ciudad: un chorro de agua que surge con una cierta inclinaci\'{o}n  describe una  \emph{par\'{a}bola}, una cadena o un cable que cuelga de dos  puntos colocados a la misma altura dibuja una \emph{catenaria},
el reflector sobre el borde de una rueda de bicicleta describe una \emph{cicloide}. El espir\'{o}grafo, juguete dise\~{n}ado por D. Fisher, que data de 1965, no traza espirales sino hermosas \emph{trocoides} (por lo que, m\'{a}s bien,  deber\'{\i}a llamarse \emph{trocoid\'{o}grafo}).
En  la ciencia f\'{\i}sica   nos topamos  con m\'{a}s curvas:  la \emph{braquist\'{o}crona} (o curva de descenso m\'{a}s r\'{a}pido), la \emph{taut\'{o}crona} (o curva en la que el tiempo de ca\'{\i}da a su punto m\'{a}s bajo no depende del punto inicial), las curvas de los  problemas  predador--presa, etc. Tambi\'{e}n abundan las curvas en las ciencias sociales. En este caso se trata de las gr\'{a}ficas de las distribuciones de las variables aleatorias, siendo la \emph{campana de Gauss}  la  m\'{a}s habitual. 
No nos debe  extra\~{n}ar la ubicuidad de las curvas pues, como dec\'{\i}a Galileo,  \lq\lq el Universo est\'{a} escrito en  lenguaje matem\'{a}tico, siendo las letras tri\'{a}ngulos, circunferencias y otras figuras geom\'{e}tricas $\ldots$\rq\rq  

\begin{figure}[ht]
 \centering
  \includegraphics[width=5cm,keepaspectratio, angle=-90]{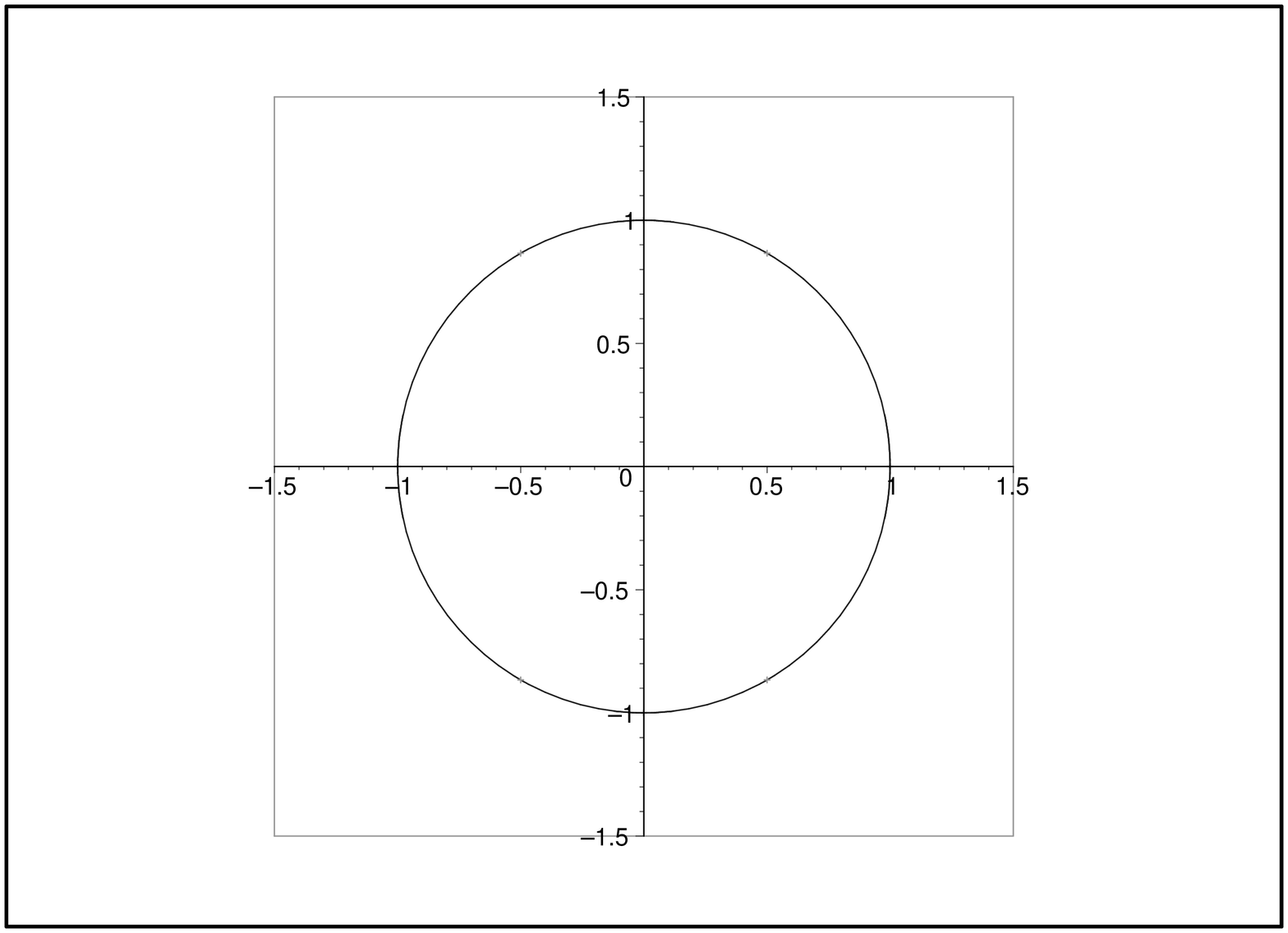}
  \caption{Circunferencia $x^2+y^2-1=0$.}
  \label{fig:fig_01}

  \includegraphics[width=5cm,keepaspectratio, angle=-90]{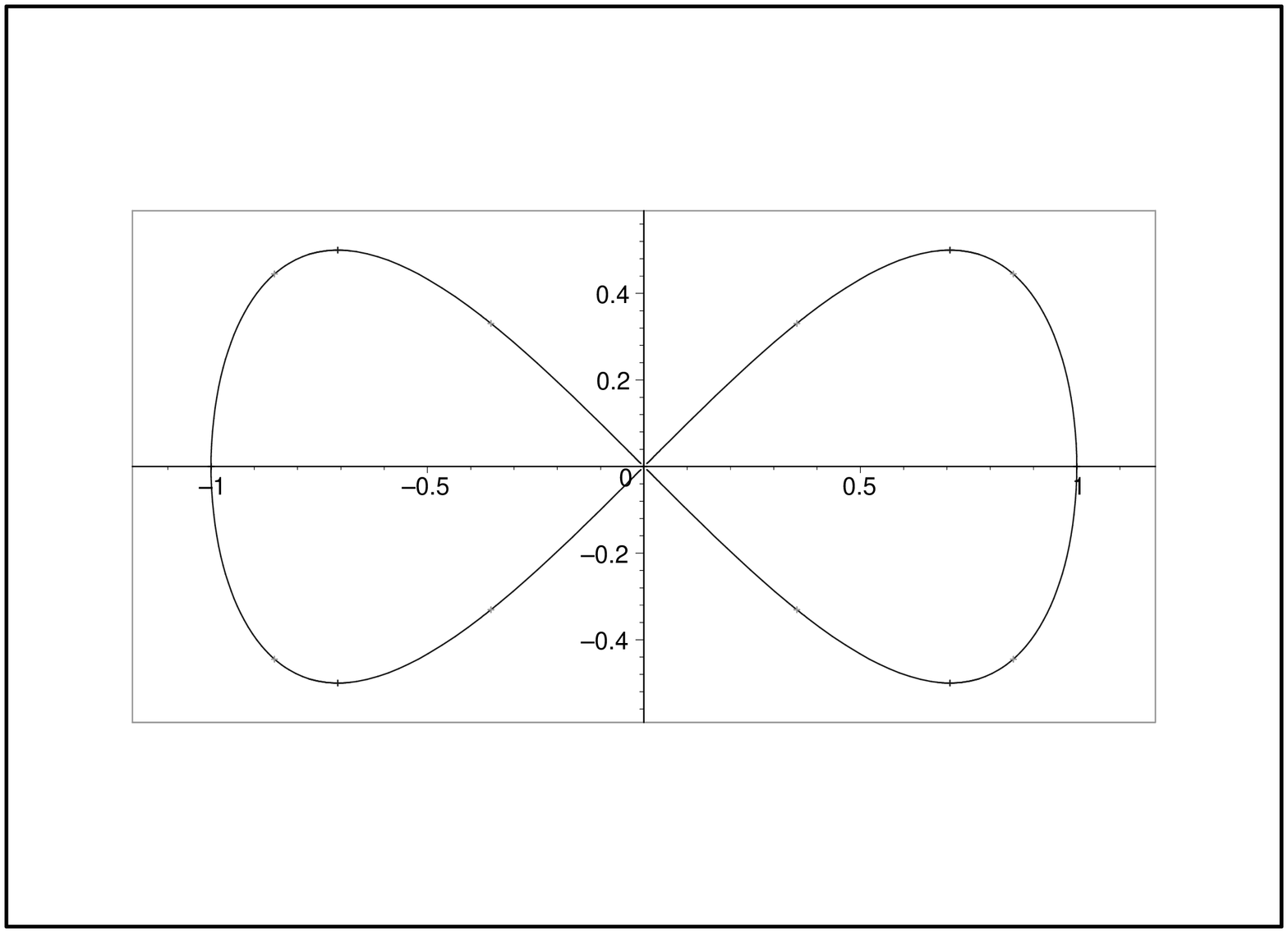}\\
  \caption{Lemniscata de Huygens.}
  \label{fig:fig_02}
  \end{figure}

El atractivo de ciertas  curvas hace que sean usadas como
 marcas, en publicidad: tres elipses tangentes en Toyota y una en Ford, Kia, Hyundai o Lexus, una par\'{a}bola en Thyssen o tres circunferencias secantes en Krupp. Por otro lado, hacemos  uso pr\'{a}ctico de algunas curvas,  por sus propiedades.  As\'{\i},  la manguera de un  extintor de incendios o una cinta cassette est\'{a}n enrolladas en espiral y  el perfil de un tornillo, un solenoide o un cable  antiguo de tel\'{e}fono tienen forma de \emph{h\'{e}lice}.

 \begin{figure}[ht]
 \centering
  \includegraphics[width=5cm,keepaspectratio, angle=-90]{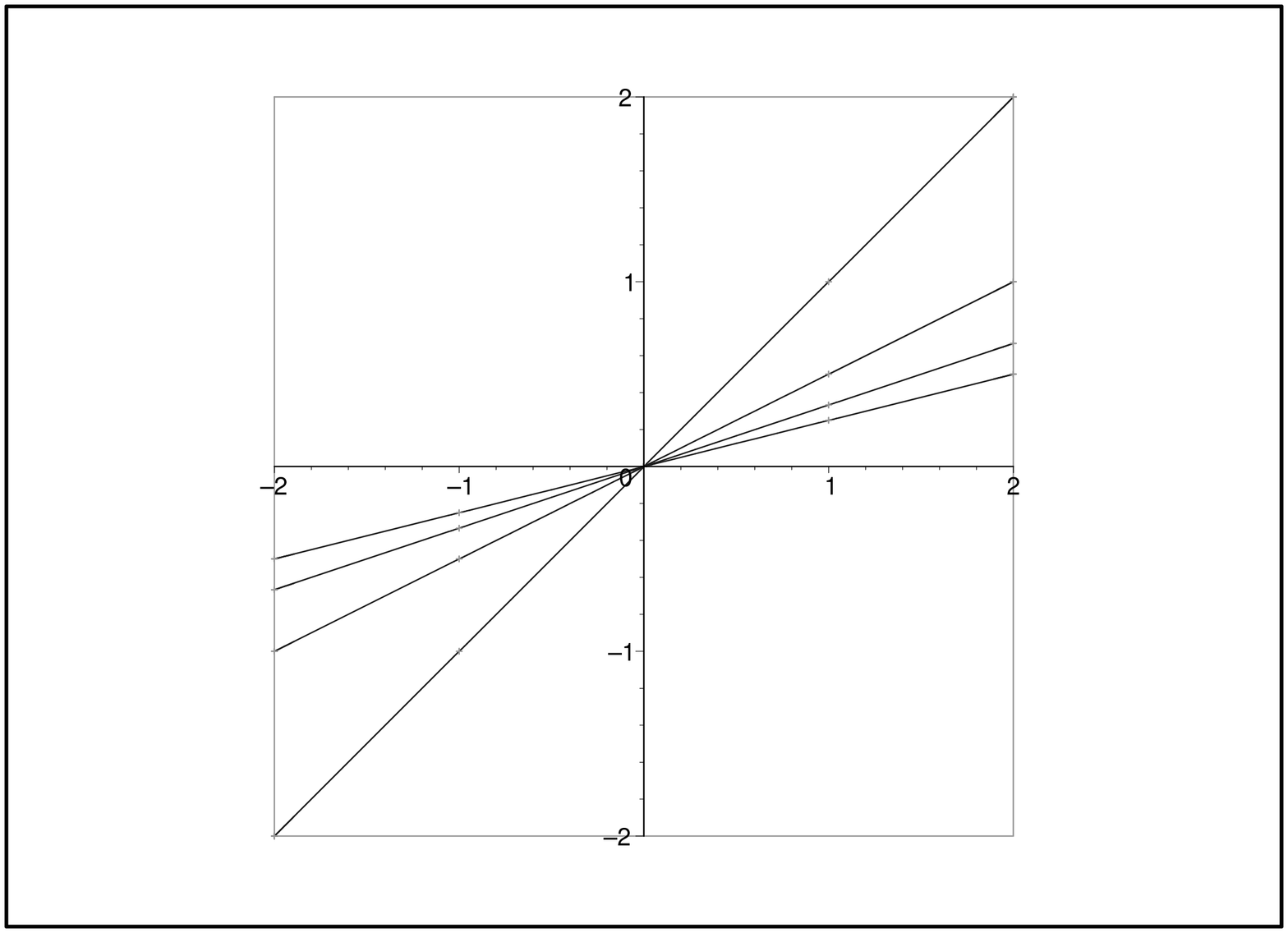}\\
  \caption{Rectas que pasan por el origen.}
  \label{fig:fig_03}

  \includegraphics[width=5cm,keepaspectratio, angle=-90]{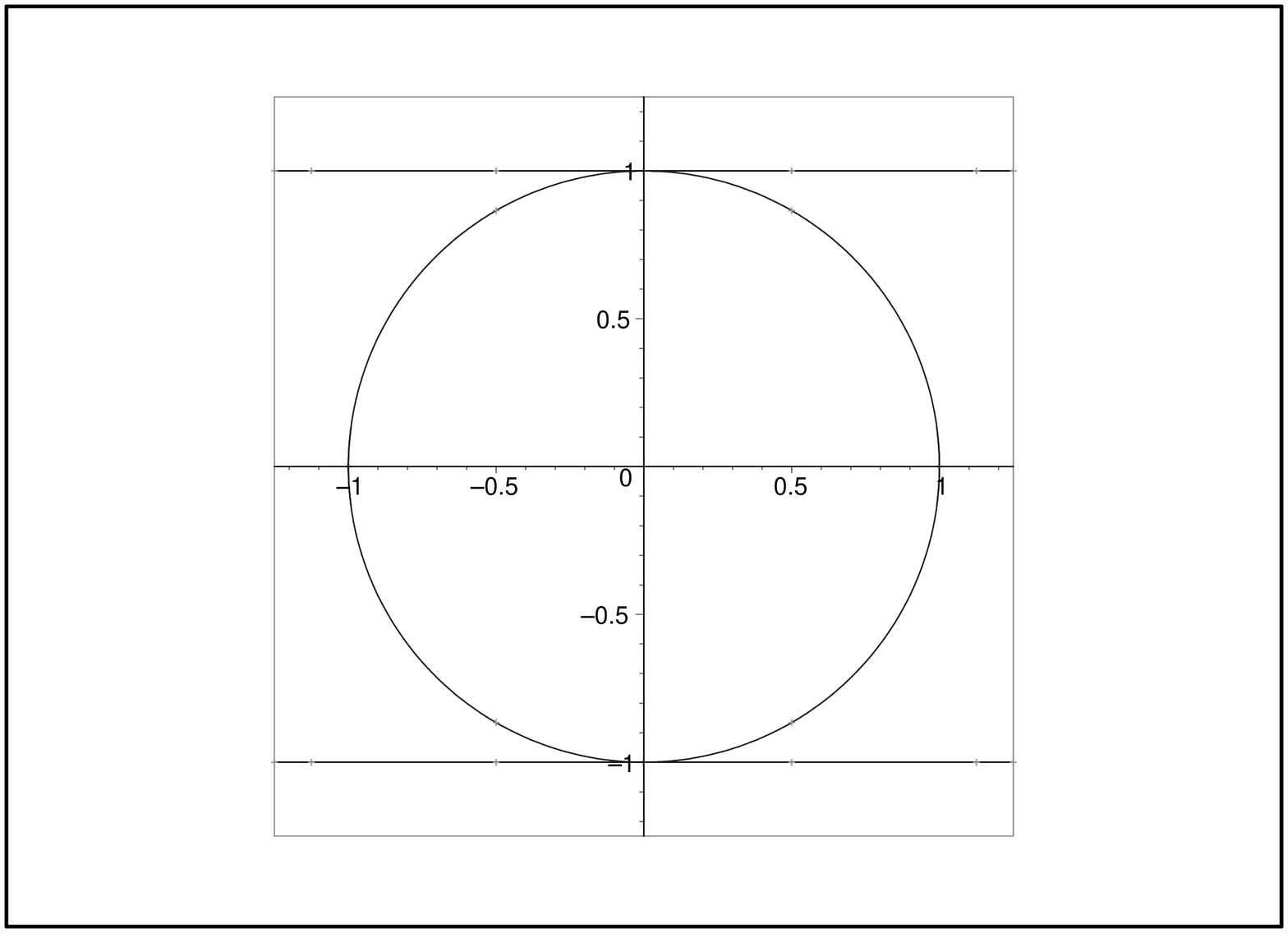}\\
  \caption{Rectas tangentes a la circunferencia en los puntos $(0,1)$ y $(0,-1)$.}
  \label{fig:fig_04}
  \end{figure}

Sin duda alguna, la \emph{circunferencia} es, tras la recta, la  curva m\'{a}s elemental. 
Bastan un punto, llamado \emph{centro}, y una cantidad positiva, llamada \emph{radio}, para describirla. Sabemos que la circunferencia $\CC$ de centro $O$ y radio uno 
es el conjunto de los puntos $P$ del plano que distan de  $O$ una unidad.

  \begin{figure}[ht]
 \centering
  \includegraphics[width=5cm,keepaspectratio, angle=-90]{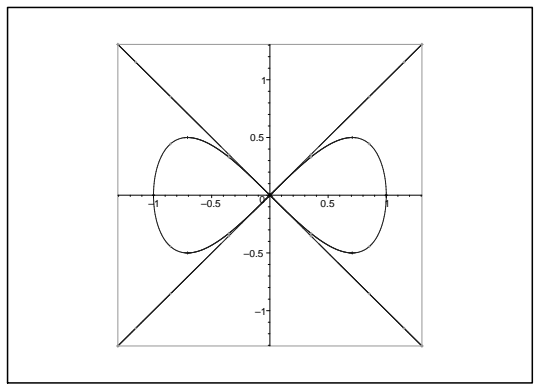}\\
  \caption{Rectas tangentes a la lemniscata en el origen.}
  \label{fig:fig_05}
  \end{figure}

  \begin{figure}
 \centering
  \includegraphics[width=5cm,keepaspectratio, angle=-90]{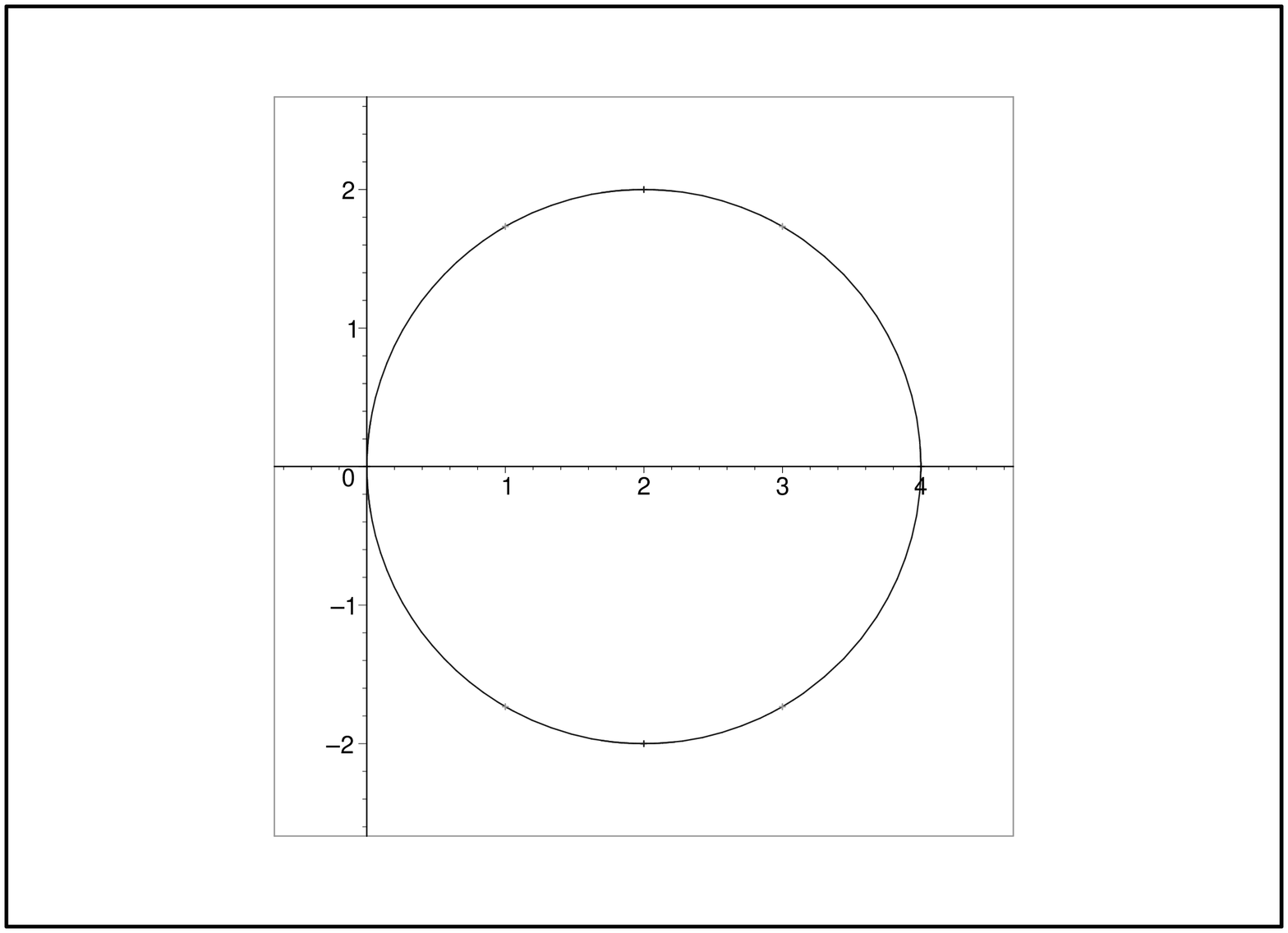}\\
  \caption{Circunferencia $x^2-4x+y^2=0$.}
  \label{fig:fig_06}
  \end{figure}


Si dotamos al plano de un sistema cartesiano 
de coordenadas 
 y  hacemos que el centro $O$ coincida con  el origen de coordenadas $(0,0)$,  entonces $\CC$ ser\'{a} el conjunto de los puntos $P(x,y)$ que  satisfacen la igualdad
$\sqrt{x^2+y^2}=1,$  ver figura \ref{fig:fig_01}. Elevando ambos miembros de la igualdad al cuadrado y pasando todos los t\'{e}rminos al primer miembro, obtenemos una condici\'{o}n equivalente
\begin{equation}\label{eqn:circunferencia}
x^2+y^2-1=0.
\end{equation}

\begin{figure}[ht]
 \centering
  \includegraphics[width=5cm,keepaspectratio, angle=-90]{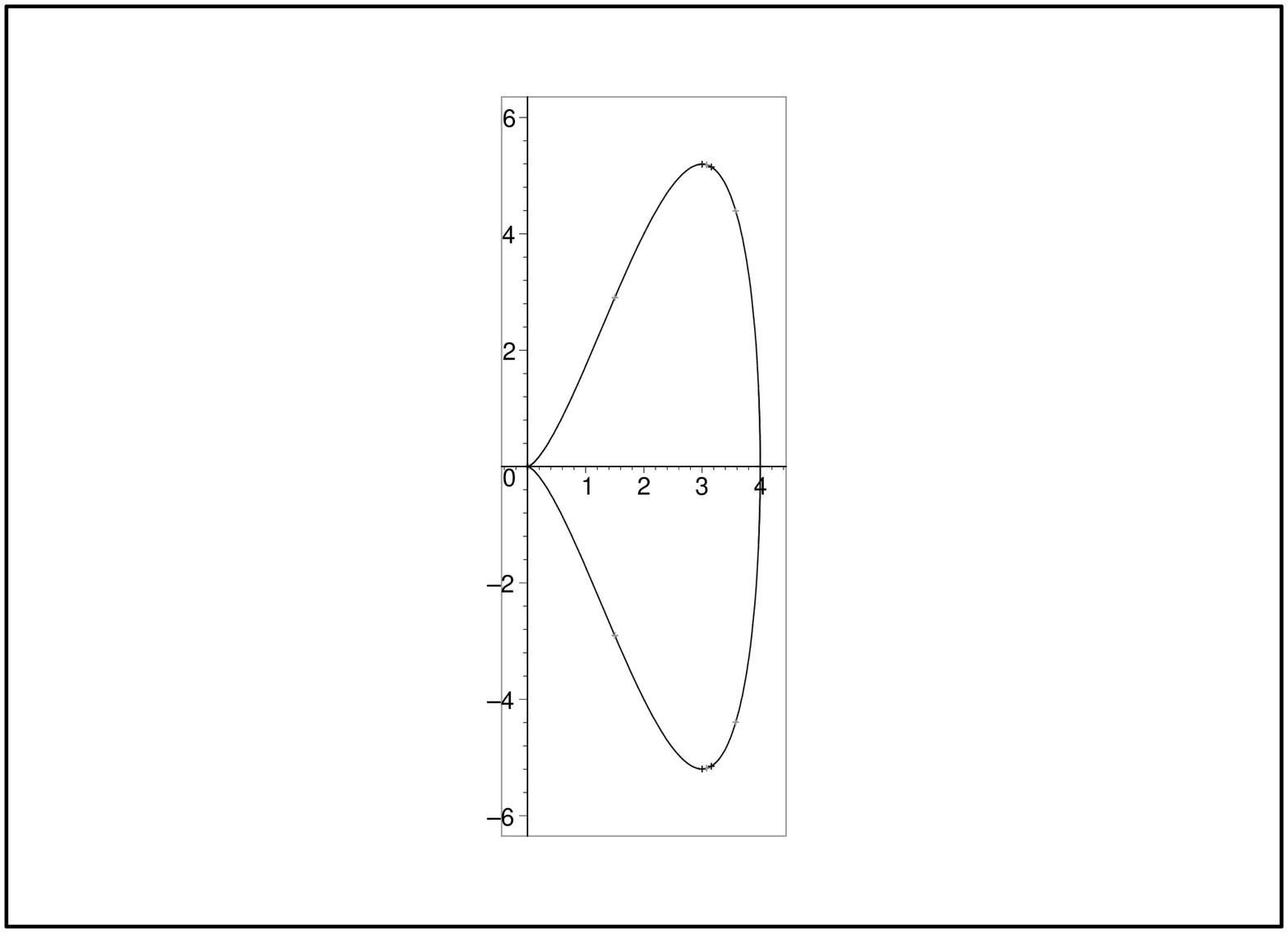}\\
  \caption{Curva piriforme.}
  \label{fig:fig_07}

  \includegraphics[width=5cm,keepaspectratio, angle=-90]{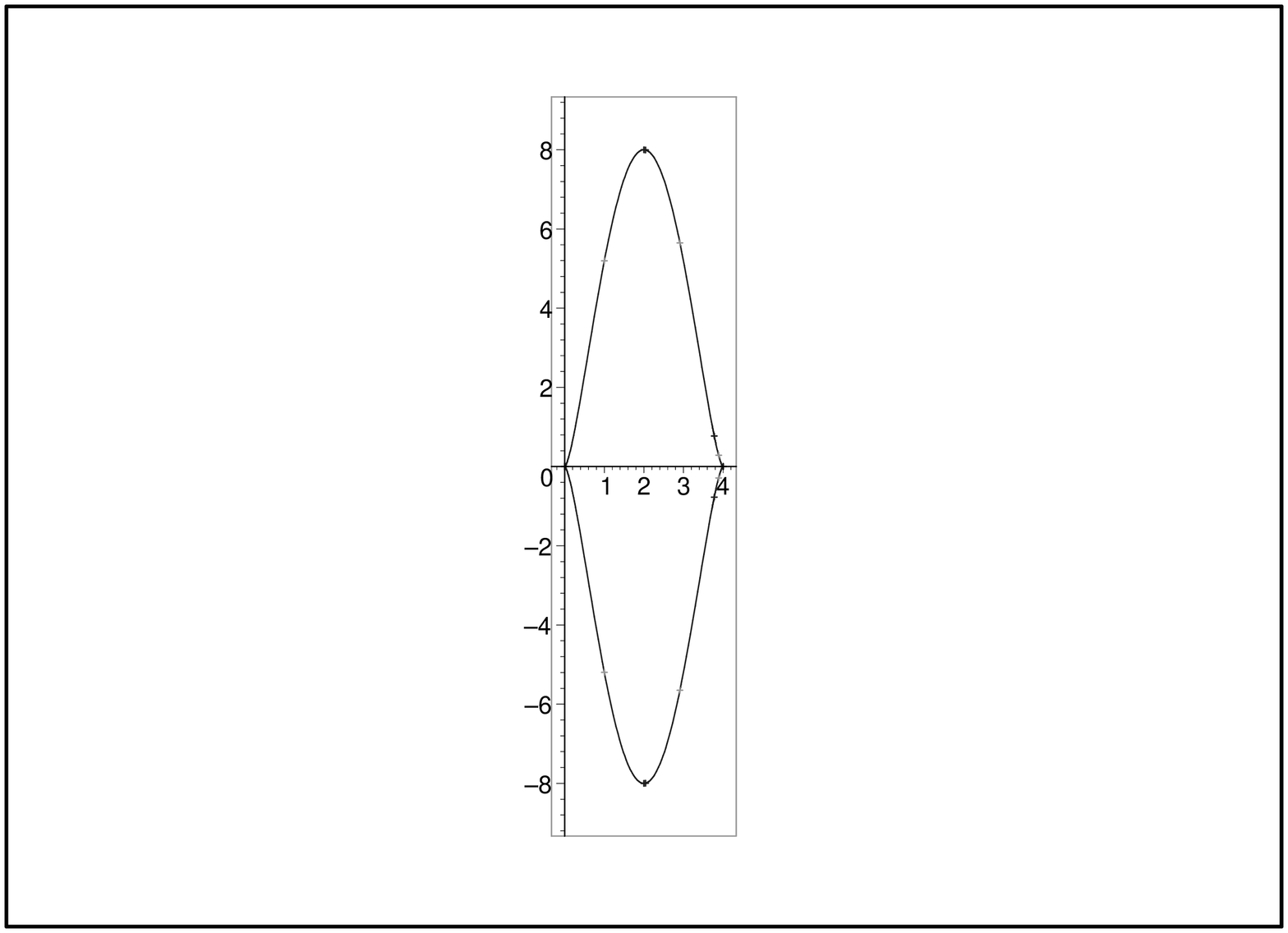}\\
  \caption{Curva en forma de labios.}
  \label{fig:fig_08}
  \end{figure}

Como las coordenadas $(x,y)$ de los puntos de $\CC$  satisfacen la \emph{relaci\'{o}n algebraica} (\ref{eqn:circunferencia}),  diremos que la circunferencia es una \emph{curva algebraica}.
No debemos creer que todas las curvas mencionadas m\'{a}s arriba son algebraicas: por ejemplo, la cicloide, la catenaria, las espirales y algunas trocoides no lo son. Tampoco son algebraicas la braquist\'{o}crona ni la taut\'{o}crona, que resultan ser cicloides invertidas.

\begin{figure}[ht]
 \centering
  \includegraphics[width=5cm,keepaspectratio, angle=-90]{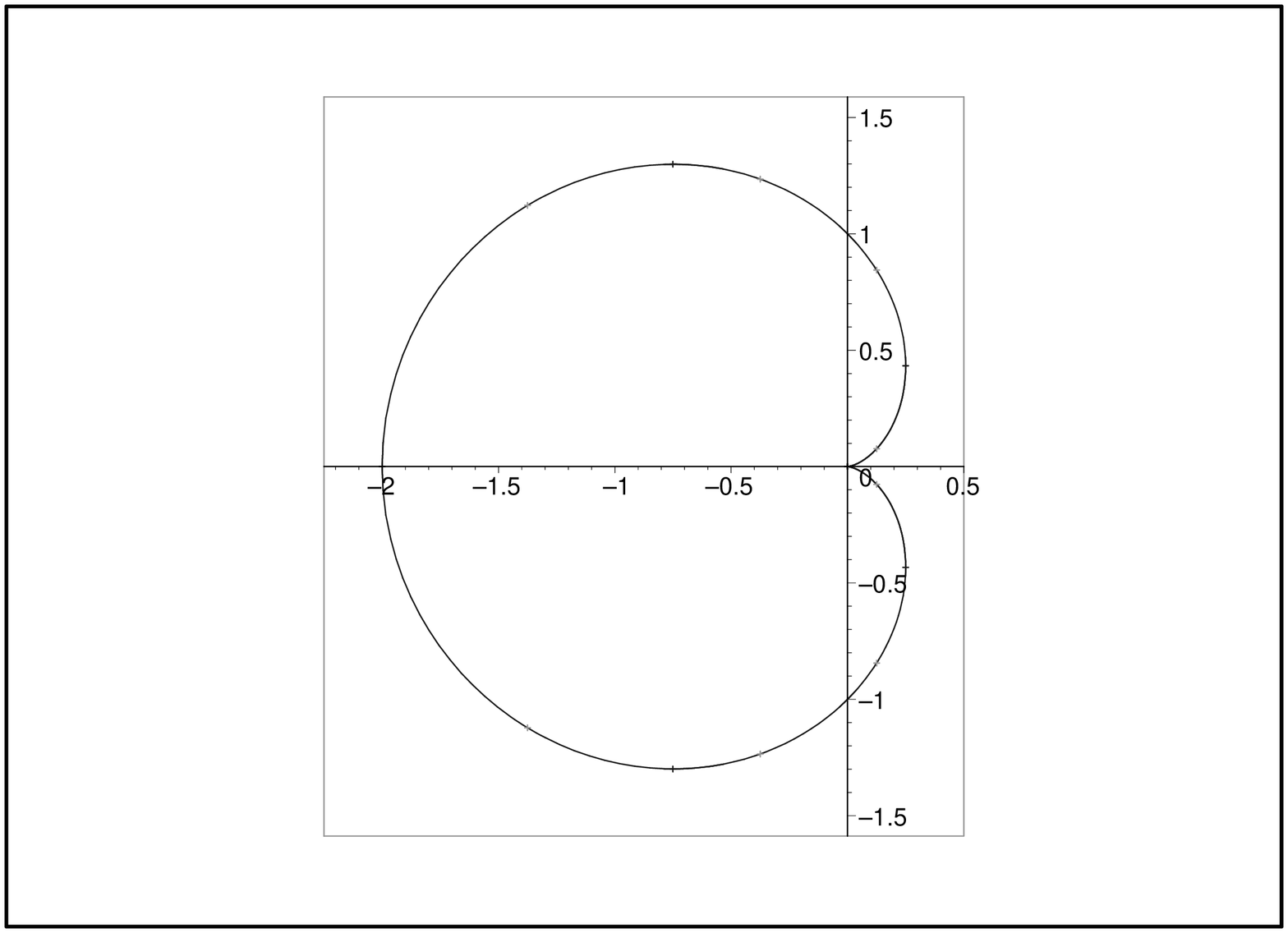}\\
  \caption{Cardioide.}
  \label{fig:fig_09}
%
  \includegraphics[width=5cm,keepaspectratio, angle=-90]{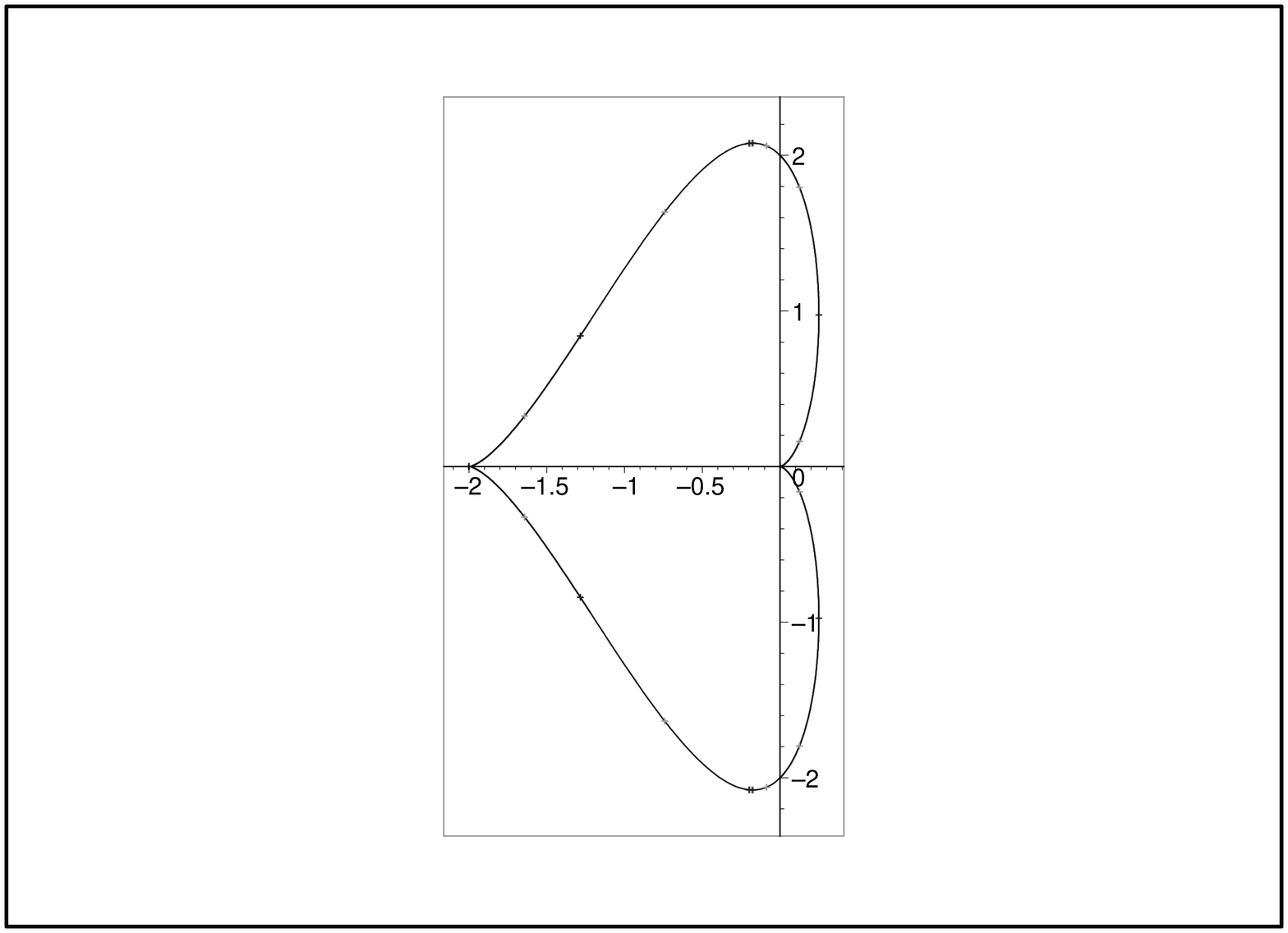}\\
  \caption{Curva en forma de coraz\'{o}n.}
  \label{fig:fig_10}
  \end{figure}

A la elegante belleza de las curvas algebraicas se suma la simplicidad de su presentaci\'{o}n: un polinomio en dos variables ($x,y$ en el caso anterior), es todo lo que necesitamos para describir una tal curva. En otras palabras, en el polinomio (sumas y productos de potencias de $x$ e $y$) est\'{a} resumida toda la geometr\'{\i}a de la curva.

\begin{figure}[ht]
 \centering
  \includegraphics[width=5cm,keepaspectratio, angle=-90]{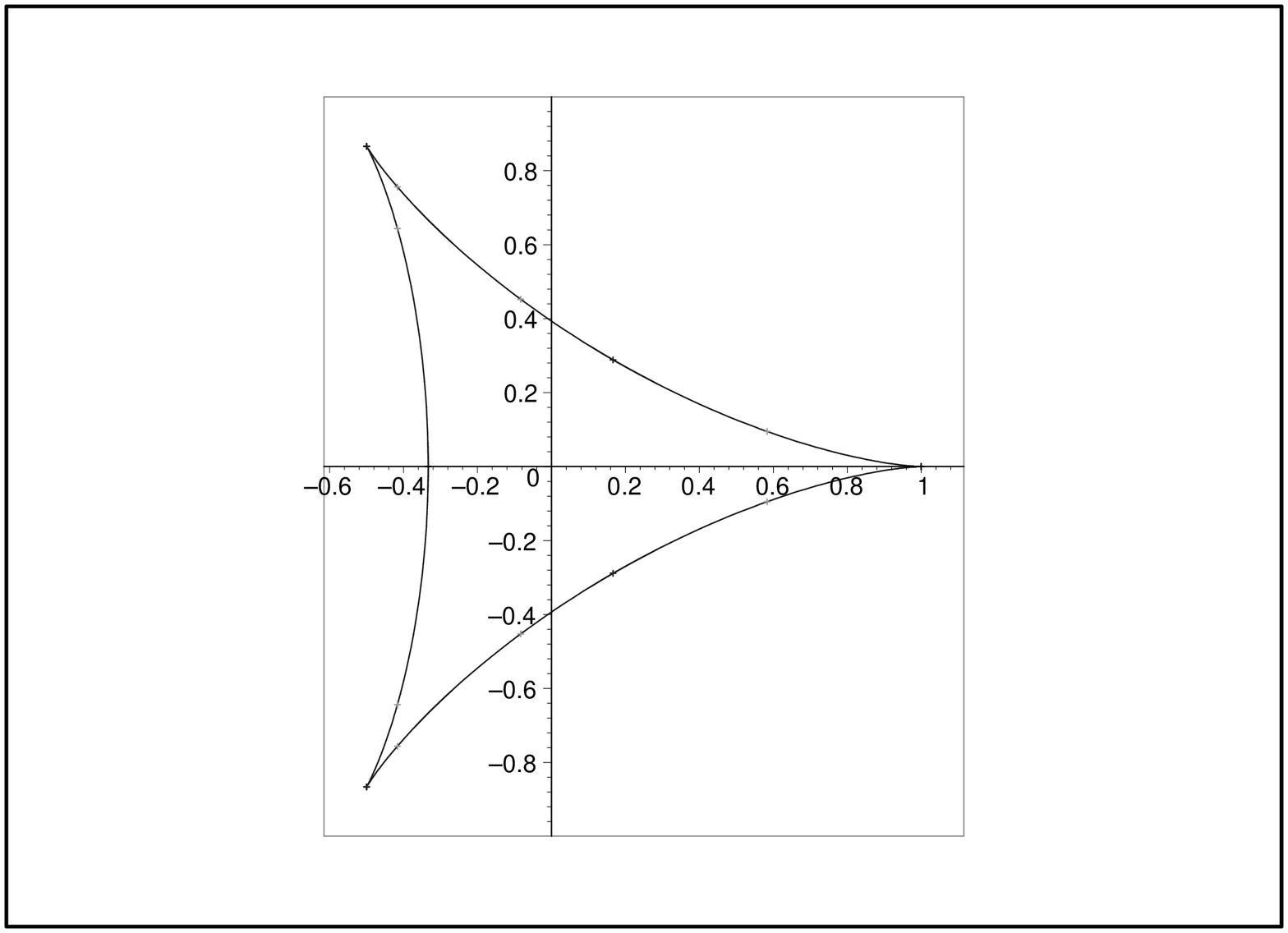}\\
  \caption{Curva tric\'{u}spide.}
  \label{fig:fig_11}

  \includegraphics[width=5cm,keepaspectratio, angle=-90]{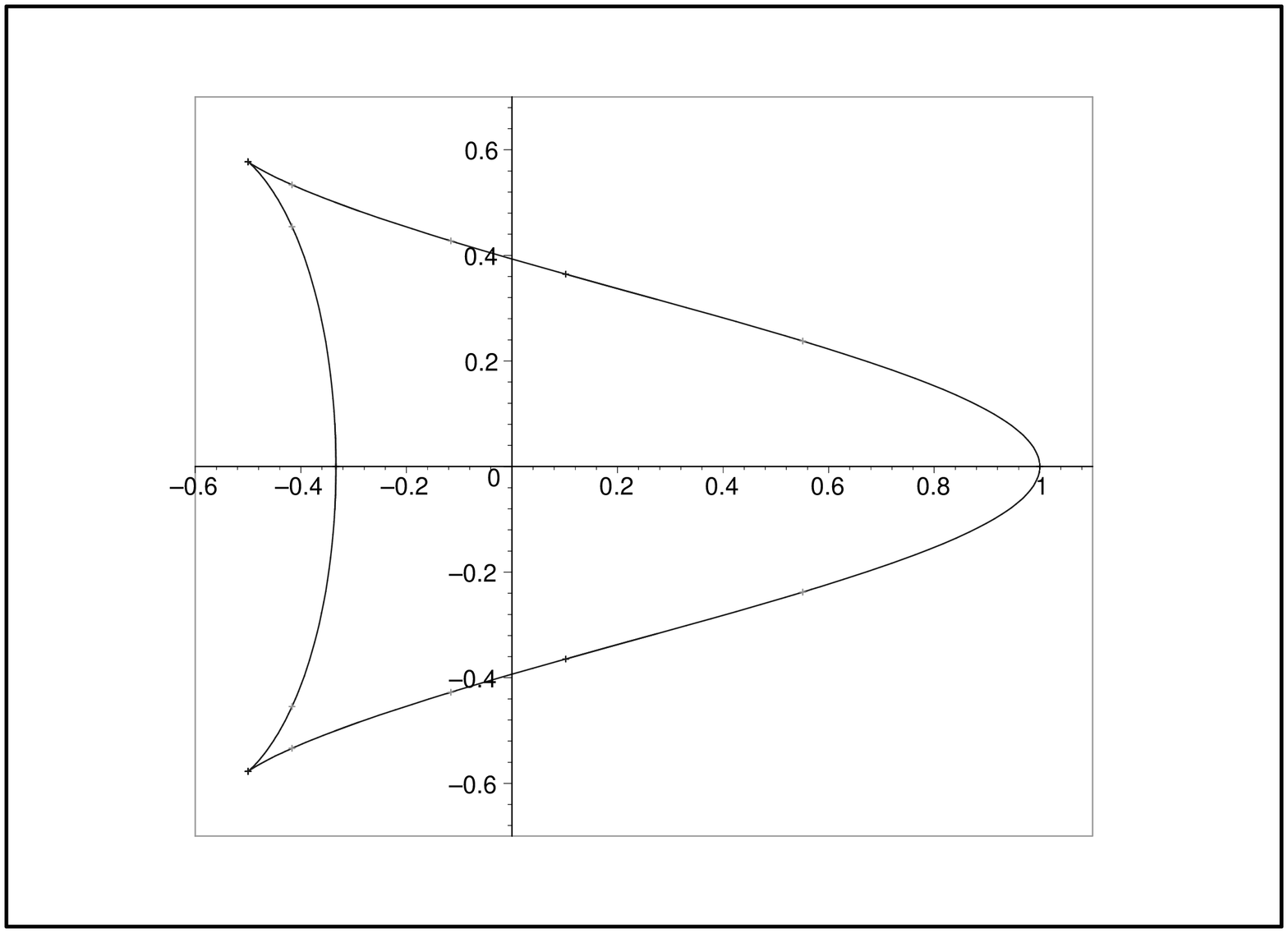}\\
  \caption{Curva en forma de punta de flecha.}
  \label{fig:fig_12}
  \end{figure}

El programa de ordenador \emph{MAPLE},  paquete \emph{algcurves}, dibuja una curva algebraica con gran precisi\'{o}n,  a  partir, exclusivamente, de su polinomio. Con \'{e}l hemos producido  las figuras de estas notas.

  \begin{figure}[ht]
 \centering
  \includegraphics[width=5cm,keepaspectratio, angle=-90]{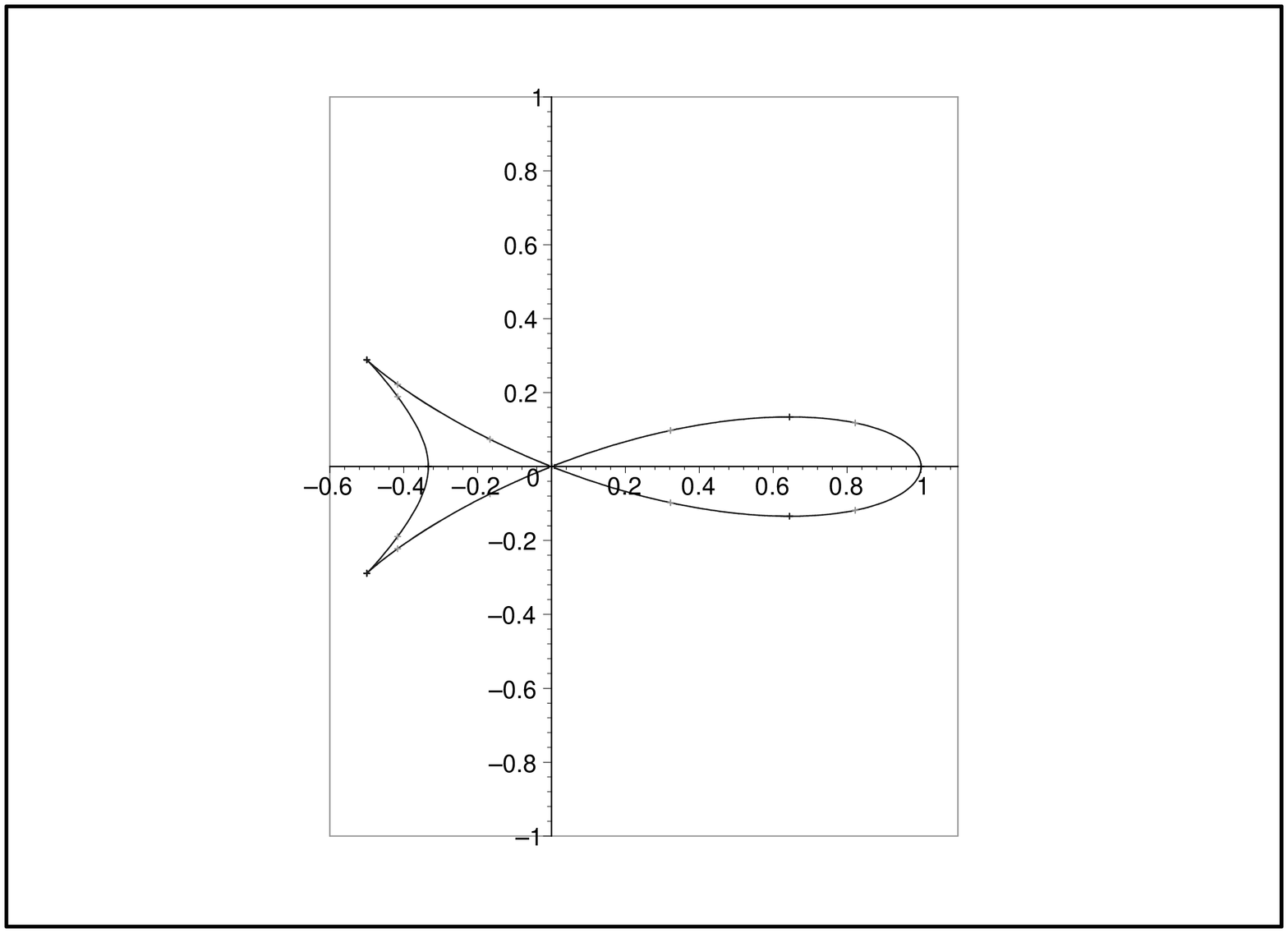}\\
  \caption{Curva pisciforme.}
  \label{fig:fig_13}
  \end{figure}

\section{Implosi\'{o}n y explosi\'{o}n}
A continuaci\'{o}n  vamos a   \emph{obtener  curvas algebraicas, con formas predeterminadas,  a partir de  circunferencias u otras curvas conocidas,} mediante dos procesos  algebraicos: \emph{la implosi\'{o}n y la explosi\'{o}n}. Con ello obtendremos, en primer lugar, algunas curvas famosas (la \emph{lemniscata de Huygens} y la curva \emph{piriforme}), para posteriormente obtener  curvas nuevas. En efecto, las curvas  con formas de \emph{labios, coraz\'{o}n, punta de flecha} y \emph{pisciforme}  se presentan  en estas notas  por vez primera, esto es, hasta donde sabemos, dichas curvas no aparecen en los tratados sobre curvas al uso.

Comencemos con la circunferencia (1). Tomamos  una variable nueva, $z$, y sustituimos $y$ por $z/x$ en la ecuaci\'{o}n (\ref{eqn:circunferencia}), obteniendo
$x^2+\left(\frac{z}{x}\right)^2-1=0.$
Multiplicando por $x^2$, para quitar denominadores,   llegamos a la ecuaci\'{o}n
\begin{equation}\label{eqn:leminscata_huy}
x^4+z^2-x^2=0.
\end{equation}

La curva $\LL$ asociada a la ecuaci\'{o}n (\ref{eqn:leminscata_huy})  se conoce como \emph{lemniscata de Huygens}, ver figura \ref{fig:fig_02}.  \emph{Cristiaan Huygens} (1629--1695) fue un matem\'{a}tico holand\'{e}s entre cuyos logros est\'{a} la patente del primer reloj de p\'{e}ndulo. Es f\'{a}cil  formar una lemniscata (que significa curva en  forma de ocho) con las manos.
Tomamos una cinta circular, sosteni\'{e}ndola con cuatro dedos de una mano y, con el pulgar y el \'{\i}ndice de  la mano libre,  pegamos dos puntos de la cinta.  Este  pegado se traduce al lenguaje matem\'{a}tico en la sustituci\'{o}n de  $y$ por  $z/x$, proceso conocido como \emph{implosi\'{o}n} (\emph{blow--down}, en ingl\'{e}s). Si $y=z/x$ entonces $z=yx$, lo que,  en el plano coordenado $XZ$,  proporciona las ecuaciones de \emph{todas} las rectas que pasan por el origen: cada recta corresponde a un valor de $y$ (que es la pendiente de dicha recta),  ver figura \ref{fig:fig_03}.

Tras la implosi\'{o}n,  los puntos $(0,1)$ y $(0,-1)$ de $\CC$ vienen a coincidir en el punto $(0,0)$ de $\LL$.
Asimismo, la recta $y=1$,  tangente a $\CC$ en el punto $(0,1)$,  se transforma en la recta $z=x$, mientras que la recta $y=-1$, va  a parar a la recta $z=-x$, ver figuras \ref{fig:fig_04} y \ref{fig:fig_05}. No es extra\~{n}o, pues  $z=x$ y $z=-x$ son las rectas tangentes a $\LL$ en el origen $(0,0)$.

El proceso contrario a la implosi\'{o}n se denomina \emph{explosi\'{o}n} (\emph{blow--up}). Se trata de sustituir $z$ por el producto $xy$, pasando  pues del plano  $XZ$ al plano $XY$. Si efectuamos el cambio  $z=xy$ en la ecuaci\'{o}n (\ref{eqn:leminscata_huy}) obtenemos
$x^4+(xy)^2-x^2=0,$ o equivalentemente,
a\begin{equation}
x^2(x^2+y^2-1)=0.
\end{equation}

 Que el  producto de $x^2$ por $x^2+y^2-1$ sea nulo, significa que  uno de los dos factores debe anularse. Por consiguiente, cada punto $(x,z)$ de $\LL$  con $x$ distinto de cero (esto es, cada  punto de $\LL$, salvo el origen) se transforma en un punto de la circunferencia  $\CC$. Adem\'{a}s,  el punto $(0,0)$ de $\LL$ \emph{explota}, convirti\'{e}ndose  en dos puntos de  $\CC$, que son $(0,1)$ y $(0,-1)$. Mediante esta explosi\'{o}n,  las rectas $z=x$ y $z=-x$, que se cortan en el origen,  se transforman  en las rectas $y=1$ e $y=-1$, que no se cortan. Este modo de \emph{eliminar un punto de auto--intersecci\'{o}n en una curva plana}  es habitual  en nuestras ciudades: mediante la construcci\'{o}n de  pasos subterr\'{a}neos, los arquitectos municipales eliminan las intersecciones de v\'{\i}as, para conseguir un tr\'{a}fico m\'{a}s fluido.
En resumen: hemos \emph{explotado} el  origen de coordenadas del plano $XZ$ y, rec\'{\i}procamente, hemos \emph{implosionado} los puntos $(0,1)$ y $(0,-1)$ del plano $XY$.

Consideremos ahora la circunferencia $\DD$ de centro $(2,0)$ y radio 2, ver figura \ref{fig:fig_06}.
Sus puntos $(x,y)$  satisfacen la igualdad
$\sqrt{(x-2)^2+(y-0)^2}=2.$ Elevando al cuadrado ambos miembros,  desarrollando los cuadrados y  pasando  todo al primer miembro, obtenemos la siguiente  ecuaci\'{o}n de $\DD$
\begin{equation}\label{eqn:circunferencia2}
x^2-4x+y^2=0.
\end{equation}

El origen  pertenece a $\DD$,  ya que  $(0,0)$ satisface la ecuaci\'{o}n (\ref{eqn:circunferencia2}). ?`En qu\'{e} se transformar\'{a} $\DD$ mediante la implosi\'{o}n  del  origen? Veamos: en (\ref{eqn:circunferencia2}), sustituimos  $y$ por $z/x$ y multiplicamos ambos miembros de la igualdad por $x^2$, para quitar denominadores, obteniendo
\begin{equation}\label{eqn:piriforme}
x^4-4x^3+z^2=0.
\end{equation}

 Se trata de  una curva \emph{piriforme} $\PP$ (curva en forma de pera, aunque, a juzgar por la gr\'{a}fica, se parece m\'{a}s a un queso gallego), ver figura \ref{fig:fig_07}.  La piriforme fue estudiada por vez primera por J. Wallis en 1685.
Como resultado de la implosi\'{o}n, en $\PP$ ha aparecido  un punto  de curioso aspecto \emph{cuspidal}, algo as\'{\i} como la comisura de unos labios. Si ahora implosionamos el punto $(4,0)$ de $\PP$, ?`conseguiremos el contorno de unos labios? Veamos: en  (\ref{eqn:piriforme}), sustituimos $z$ por $t/(x-4)$, obteniendo
$x^4-4x^3+\left(\frac{t}{x-4}\right)^2=0.$ Quitando denominadores, operando y simplificando, llegamos a la ecuaci\'{o}n
\begin{equation}
x^6-12x^5+48x^4-64x^3+t^2=0,
\end{equation}
que describe una curva algebraica del plano $XT$ con la forma de labios (pronunciando la u francesa) buscada, ver figura \ref{fig:fig_08}.

  La curva  \emph{cardioide} (en forma de coraz\'{o}n) fue considerada por  O.C. Roemer, en 1676, y P. de la Hire en 1708, ver figura \ref{fig:fig_09}.  
Se trata de la trayectoria de un punto colocado sobre una circunferencia $\EE$  que rueda sin deslizamiento alrededor de otra circunferencia fija $\EE'$,  de id\'{e}ntico radio. Es f\'{a}cil  visualizar esta curva, usando dos monedas de igual valor.  En el borde de una de ellas  marcamos un punto. Dejando  fija la moneda no pintada, hacemos girar una sobre otra y el punto  describir\'{a} una \emph{cardioide} $\CC$. Su ecuaci\'{o}n es
\begin{equation}\label{eqn:cardioide_cartes}
(x^2+y^2+x)^2-x^2-y^2=0.
\end{equation}
A pesar de su nombre, $\CC$ no se parece mucho a la silueta  de  coraz\'{o}n del imaginario colectivo, pues el punto  $(-2,0)$, que pertenece a $\CC$, deber\'{\i}a ser  cuspidal. ?`Conseguiremos la forma de coraz\'{o}n deseada efectuando  una implosi\'{o}n del punto $(-2,0)$? Veamos: sustituimos  $y$ por $z/(x+2)$, obteniendo
$\left(x^2+\left(\frac{z}{x+2}\right)^2+x\right)^2-x^2-\left(\frac{z}{x+2}\right)^2=0,$ donde,  
operando y simplificando llegamos a la larga expresi\'{o}n siguiente, cuya curva est\'{a} representada en la figura   \ref{fig:fig_10},
\begin{equation}\label{eqn:corazon}\begin{split}
x^8+10x^7+40x^6+80x^5+2x^4z^2+80x^4+32x^3\\+10x^3z^2+15x^2z^2+4xz^2+z^4-4z^2=0.
\end{split}\end{equation}

La curva \emph{hipocicloide de tres puntas} o \emph{deltoide} (por su parecido con la letra griega delta may\'{u}scula) fue estudiada en profundidad por L. Euler (1707--1783) y J. Steiner  (1796--1863).   Es la
trayectoria de un punto colocado sobre una circunferencia $\EE$  que rueda sin deslizamiento por el interior de otra circunferencia fija $\EE'$,  cuando el radio de $\EE$ es la tercera parte del de $\EE'$. Podemos ver esta curva  usando  dos monedas, una de las cuales debe tener un radio tres veces mayor que la otra. Tambi\'{e}n se denomina \emph{tric\'{u}spide}, por poseer tres puntos cuspidales, ver figura \ref{fig:fig_11}.
Su ecuaci\'{o}n  es
\begin{equation}\label{eqn:tricuspide}
3(x^2+y^2)^2+8x(3y^2-x^2)+6x^2+6y^2-1=0.
\end{equation}

A partir de la tric\'{u}spide $\TT$,  vamos a encontrar una curva \emph{pisciforme}, esto es, en forma de pez. Primero explotamos en $\TT$ el punto $(1,0)$, (?`sabr\'{\i}as hacerlo, lector?)  obteniendo  una curva $\FF$ en forma de \emph{punta de flecha redondeada}, de ecuaci\'{o}n
\begin{equation}\label{eqn:punta_de_flecha}
3x^2z^4+6x^2z^2+3x^2-6xz^4+24xz^2\\-2x+3z^4+6z^2-1=0,
\end{equation}
ver figura \ref{fig:fig_12}.
A continuaci\'{o}n, implosionamos dos puntos de esta curva, haciendo $z=t/x$ y obtenemos  la ecuaci\'{o}n
\begin{equation}\label{eqn:pisciforme}
3x^6-2x^5+6x^4t^2-x^4+24x^3t^2\\+3x^2t^4+6x^2t^2-6xt^4+3t^4=0,
\end{equation}
cuya  curva es la \emph{pisciforme} buscada, ver figura \ref{fig:fig_13}.

\section{Conclusiones}
Ahora dispones, lector, de una doble herramienta (la implosi\'{o}n y la explosi\'{o}n) para crear nuevas curvas a partir de curvas conocidas, seg\'{u}n tu propio gusto.


\section*{Agradecimientos} Este art\'{\i}culo ha sido redactado con el apoyo del proyecto investigador UCM 910444.

\nocite{*}
\bibliographystyle{plain}
\bibliography{biblio_c_p_a}
\end{document}